# Plotting positions close to the exact unbiased solution: application to the Pozzuoli's bradyseism earthquake data

Pasquale Erto and Antonio Lepore

*University of Naples Federico II, Department of Industrial Engineering, Naples, Italy*

**Summary.** Graphical techniques are recommended for critical applications in order to share information with non-statisticians, since they allow for a visual analysis and helpful understanding of the results. However, graphical estimation methods are often underestimated because of their minor efficiency with respect to the analytical ones. Therefore, finding unbiased plotting positions can contribute to rise their reputation and to encourage their strategic use. This paper proposes a new general plotting position formula which can be as close as needed to the exact unbiased plotting positions. The ability of the new solution in estimating quantiles for both symmetrical and skewed location-scale distributions is shown via Monte Carlo simulation. An applicative example shows how the proposed formula enables to perform, with known accuracy, the graphical analysis of critical data, such as the earthquake magnitudes registered during the serious 1983-1984 bradyseismic crisis in Campi Flegrei (Italy). Moreover, the proposed formula gives a unified look at existing plotting positions and a definitive insight into plotting position controversies recently renewed in the literature.

## 1. Introduction

As confirmed by the renewed interest appeared in the recent literature (Rigdon and Basu 1989, Makkonen 2006, de Haan 2007, Makkonen 2008a, Cook 2011, Cook 2012, Kim *et al.* 2012, Erto and Lepore 2013, Fuglem *et al.* 2013, Makkonen 2013, Lozano-Aguilera *et al.* 2014) practitioners are used to exploiting modern software that adopts graphical estimation methods, even if there is a variety of effective analytical methods available, such as Maximum Likelihood and Bayesian techniques. In fact, especially in critical applications, the graphical estimation gives the unique opportunity to share statistical information with non-statisticians by allowing a visual check of the fit of the chosen model and by giving helpful understanding of the consequent conclusions. Clearly, if the approach is to be purely analytical there is no point in using a probability paper (Kimbal 1960).

Plotting positions have been used and discussed for many years by engineers, hydrologists and statisticians. Noticeable remarks on classical extreme value analysis and plotting positions are included in Harris (1996), Palutikof *et al.* (1999), Filliben (2001), Folland and Anderson (2002), Cook *et al.* (2003), Rasmussen and Gautam (2003), Whalen *et al.* (2004), Cook and Harris (2004), McRobie (2004), Jordan (2005), Kharin and Zwiers (2005), Kidson and Richards (2005). A comprehensive review of the main plotting positions can be found in Harter (1984), which concludes that the distribution of the variable under consideration (i.e., the parent distribution) and the purpose for which the results are to be used are the major factors in the choice of plotting position. According to this conclusion, a unique formula that outstands in any condition does not exist. In fact, to mention but a few, MATLAB software adopts Hazen's formula (Hazen 1914) as default plotting position; the National Institute of Standards and Technology handbook of statistical methods recommends Jenkinson formula (Jenkinson 1969) as well as Gringorten formula (Gringorten 1963) is preferred by Palutikof *et al.* (1999) in calculating extreme wind speeds; Hazen and Cunnane formulas are recommended by Jordaan (2005).

As mentioned above, the issue of determining a unique and distribution-free plotting position formula has recently come to light again (Lozano-Aguilera *et al.* 2014, Erto and Lepore 2013, Makkonen 2008a, Makkonen 2008b). The classical distribution-free plotting position proposed by Meeker and Escobar (1998) and promoted by Gumbel (1958) (known as the Weibull formula) is indicated by Makkonen (2008b) as the only correct solution by concluding that "the plotting position in the extreme value analysis should be considered not as an estimate, but to be equal to $i/(N+1)$ regardless of the parent distribution and the application".

This sharp statement has been confuted by many authors de Haan (2007), Makkonen (2007), Makkonen (2011), Cook (2011), Cook (2012), Erto and Lepore (2013), Fuglem *et al.* (2013), giving rise to a wide controversial discussion. However, as detailed in Section 2, most of the argumentations addressed in such discussion were already clear to Hahn and Shapiro (1967).

## 2. The plotting position controversy

### 2.1. Distribution-free versus exact unbiased plotting positions

If the distribution of the variable under consideration is known, EUPP (Exact Unbiased Plotting Positions) can be obtained via order statistics-theory. Unfortunately, in many cases these are too complex to determine (see, e.g., Lieblein and Salzer 1957) and cannot be used for practical purposes such as the return period estimation (see conclusion 4 by Cunnane 1978 and motives 12 by Lozano-Aguilera *et al.* 2014).

The practical form

$$\hat{F}_i = \frac{i - A}{N + B} \qquad i = 1, ..., N \tag{1}$$

is widely utilized to get an approximation of the EUPP (Gringorten 1963, Cunnane 1978, Guo 1990) by choosing suitable real constants $A$ and $B$ (Table 1).

Most of the plotting positions appeared in the literature (Table 2) reduce to

$$\hat{F}_i = \frac{i - A}{N + 1 - 2A} \tag{2}$$

which is obtained from (1) upon setting $B = 1 - 2A$ (Blom 1958).

It can be easily shown that (2) implies the following assumption

$$\hat{F}_i = 1 - \hat{F}_{N-i+1}. \tag{3}$$

which, if $N$ is odd, includes the results $\hat{F}_{(N+1)/2} = 1/2$, already stated by Erto and Lepore (2013).

Distribution-free approaches are also noticeable in the literature: Gumbel (1958), Makkonen (2006), Erto and Lepore (2013), Lozano-Aguilera *et al.* (2014). Most of them are essentially based on the median or the mean value of the cdf $F_X(X_{(i)})$ which, apart from the parent distribution, can be shown to be a Beta random variable $U_{(i)}$ with probability density function (pdf)

$$f_{U_{(i)}}(t) = \frac{\Gamma(a+b)}{\Gamma(a)\Gamma(b)} t^{a-1}(1-t)^{b-1} \tag{4}$$

where $a = i$ and $b = N - i + 1$.

Makkonen (2008a) develops his distribution-free approach by interpreting the plotting position as the non-exceedance probability of the next observation in an order ranked sample $P\{X \leq X_{(i)}\}$ and obtains (Makkonen *et al.* 2013)

$$\hat{F}_i = P\{X \leq X_{(i)}\} = E\{F_X(X_{(i)})\} = \frac{i}{N+1} \tag{5}$$

widely known as the Weibull plotting position (Gumbel 1958).

### 2.2. Exact Unbiased Plotting Position approach

If $X$ (and then $X_{(i)}$) is a continuous location-scale random variable, we can introduce the reduced variate

$$Z_{(i)} = (X_{(i)} - a)/b \tag{6}$$

where $a$ and $b$ are the location and the non-negative scale parameters, respectively. Obviously from (4)

$$F_Z(Z_{(i)}) = F_X(X_{(i)}) = U_{(i)}. \tag{7}$$

From (6) and (7) it then follows that

$$E\{F_Z^{-1}(U_{(i)})\} = E\{Z_{(i)}\} = (E\{X_{(i)}\} - a)/b. \tag{8}$$

The formula proposed in this paper is based on the Taylor series expanding of $F_Z^{-1}(U_{(i)})$ around $E\{U_{(i)}\} = \mu_{U_{(i)}}$

$$E\{Z_{(i)}\} = \sum_{j=0}^{\infty} \frac{F_Z^{-1}(\mu_{U_{(i)}})}{j!} E\left\{(U_{(i)} - \mu_{U_{(i)}})^j\right\}. \tag{9}$$

In order to graphically estimate the distribution (location and scale) parameters through probability papers the following regression model is assumed from (8)

$$x_{(i)} = b y_{(i)} + a + \varepsilon_{(i)} \tag{10}$$

where $\varepsilon_{(i)}$ represent the error/residual.

Apart from the method used to estimate $a$ and $b$ through probability paper, the plotting positions proposed in the past decades are nothing but different formulas used to obtain approximations for $y_{(i)}$ (generally in the practical form (1) or (2)).

In accordance with the Cunnane (1978) proposal, if we assume $y_{(i)} = E\left\{Z_{(i)}\right\}$, then the covariance $\sigma_{(X_{(i)}, X_{(j)})}$ between $X_{(i)}$ and $X_{(j)}$ is nonzero and can be expressed in term of the covariance $\sigma_{(i,j)}$ between $Z_{(i)}$ and $Z_{(j)}$ as follows

$$\sigma_{(X_{(i)}, X_{(j)})} = b^2 \sigma_{(i,j)}. \tag{11}$$

Therefore, the covariance matrix of the error $\boldsymbol{\varepsilon} = \left[\varepsilon_{(1)} \cdots \varepsilon_{(N)}\right]'$

$$\mathbf{V} = \begin{bmatrix} \sigma^2_{(1)} & \dots & \sigma_{(1,N)} \\ \vdots & \sigma_{(i,j)} & \vdots \\ \sigma_{(N,1)} & \dots & \sigma^2_{(N)} \end{bmatrix} \tag{12}$$

where $\sigma^2_{(i)} = \sigma_{(i,i)}$, has nonzero off-diagonal elements and different diagonal elements and can be shown to be nonsingular and positive definite.

In matrix notation, the regression model can be then expressed as

$$\mathbf{X} = \mathbf{A}\boldsymbol{\theta} + \boldsymbol{\varepsilon} \tag{13}$$

where $\boldsymbol{\theta} = (a, b)$ and the $n \times 2$ matrix

$$\mathbf{A} = \begin{bmatrix} 1 & E\left\{Z_{(1)}\right\} \\ \vdots & \vdots \\ 1 & E\left\{Z_{(N)}\right\} \end{bmatrix}. \tag{14}$$

The generalized least-squares solution to the regression model (13) is

$$\hat{\boldsymbol{\theta}} = \begin{bmatrix} \hat{a} \\ \hat{b} \end{bmatrix} = \left(\mathbf{A}'\mathbf{V}^{-1}\mathbf{A}\right)^{-1}\mathbf{A}'\mathbf{V}^{-1}\mathbf{X} \tag{15}$$

which can be shown to be the Best Linear Unbiased Estimators (BLUEs) of $\boldsymbol{\theta}$ (Lieblein 1953, Draper and Smith 1981). In the case of the Uniform distribution, these estimators are also Minimum Variance Unbiased Estimator (MVUE).

Now it is clear in which sense the Cunnane proposal is the EUPP and can be obtained as shown by Cunnane (1978) and recently encouraged by Hong and Li (2013) and Fuglem *et al.* (2013).

However, not to miss the applicative use of probability papers, Hong and Li (2013) as well Hahn and Shapiro (1967) remark that if the sample size is not dramatically small the Ordinary Least Squares (OLS) method (which is the nearest method to the visual best fitting track of a line) gives satisfactory estimates.

The simulation study proposed in the next section shows that, for practical purposes, satisfactory results can be achieved by the solution $\hat{F}_i = F_Z\left(y_{(i)}\right)$ obtained by considering the first $k = 4$ terms of (9) in the model (10)

$$y_{(i)} \simeq F_Z^{-1}\left(\mu_{U_{(i)}}\right) + \frac{\mu_{U_{(i)}}\left(1-\mu_{U_{(i)}}\right)}{2\left(N+2\right)} F_Z^{-1(2)}\left(\mu_{U_{(i)}}\right) + \frac{\mu_{U_{(i)}}\left(1-\mu_{U_{(i)}}\right)}{\left(N+2\right)^2}$$
$$\left\{\frac{1}{3}\left(1-2\mu_{U_{(i)}}\right)F_Z^{-1(3)}\left(\mu_{U_{(i)}}\right) + \frac{1}{8}\mu_{U_{(i)}}\left(1-\mu_{U_{(i)}}\right)F_Z^{-1(4)}\left(\mu_{U_{(i)}}\right)\right\} \tag{16}$$

and

$$\sigma_{(i,j)} \simeq \frac{\mu_{U_{(i)}}\left(1-\mu_{U_{(j)}}\right)}{N+2} F_Z^{-1}\left(\mu_{U_{(i)}}\right) + \frac{\mu_{U_{(i)}}\left(1-\mu_{U_{(j)}}\right)}{\left(N+2\right)^2}\left\{\left(1-2\mu_{U_{(i)}}\right)F_Z^{-1(2)}\left(\mu_{U_{(i)}}\right)\right.$$
$$F_Z^{-1}\left(\mu_{U_{(j)}}\right) + \left(1-2\mu_{U_{(j)}}\right)F_Z^{-1(2)}\left(\mu_{U_{(j)}}\right)F_Z^{-1}\left(\mu_{U_{(i)}}\right) + \frac{1}{2}\mu_{U_{(i)}}\left(1-\mu_{U_{(i)}}\right)$$
$$F_Z^{-1(3)}\left(\mu_{U_{(i)}}\right)F_Z^{-1}\left(\mu_{U_{(j)}}\right) + \frac{1}{2}\mu_{U_{(j)}}\left(1-\mu_{U_{(j)}}\right)F_Z^{-1(3)}\left(\mu_{U_{(j)}}\right)F_Z^{-1}\left(\mu_{U_{(i)}}\right)$$
$$\left. + \frac{1}{2}\mu_{U_{(i)}}\left(1-\mu_{U_{(j)}}\right)F_Z^{-1(2)}\left(\mu_{U_{(i)}}\right)F_Z^{-1(2)}\left(\mu_{U_{(j)}}\right)\right\}. \tag{17}$$

where $F_Z^{-1(j)}\left(\cdot\right)$ is the $j$-th derivative of $F_Z^{-1}\left(\cdot\right)$. Note that if $k = 0$ and $\mathbf{V} = \sigma\mathbf{I}$, where $\mathbf{I}$ is the identity matrix of order $N$, the proposed solution in (16) coincides with the Weibull plotting position promoted by Makkonen (2008b).

**Table 1: Most relevant EUPP approximations in the form (1) or (2)**

| Author(s) | Distribution | *A* | *B* |
|---|---|---|---|
| Hazen (1914) - Foster (1936) | Gumbel | 1/2 | $1-2A$ |
| Beard (1943) | Normal | 0.31 | $1-2A$ |
| Blom (1958) | Normal | 3/8 | $1-2A$ |
| Tukey (1962) | Normal | 1/3 | $1-2A$ |
| Gringorten (1963) | Gumbel | 44 | $1-2A$ |

| Yu and Huang (1999) (a) | Normal | 0.399 | 0.203 |
| Yu and Huang (1999) (b) | Gumbel | 0.507 | 0.176 |
| De (2000) | Gumbel | 0.28 | 0.28 |

**Table 2: Most relevant distribution-free plotting positions in the form (1)**

| Author(s) | $A$ |
|---|---|
| Weibull (1914) | 0 |
| Cunnane (1977) | 2/5 |
| Adamowski (1981) | 1/4 |
| Kerman (2011) | 1/3 |
| Erto and Lepore (2013) | $N+(N-1)/(2^{1/N}-2)$ |

## 3.  Monte Carlo simulation for testing plotting position descriptive and predictive ability

The linear estimators $\hat{a}$ and $\hat{b}$ obtained from (10) for location-scale distributions (and related families) are equivariant (Erto 1981) and therefore, the quantities

$$K_1 = \frac{a-\hat{a}}{\hat{b}} \quad \text{and} \quad K_2 = \frac{\hat{b}}{b} \tag{18}$$

are parameter-free (Lawless 1978). Moreover, if we denote with $x_T$ the theoretical quantile of $X$ at the given return period $T$ and with $\hat{x}_T$ its graphical (OLS) estimate, it can be readily shown that also

$$K_3 = \frac{\hat{x}_T - x_T}{b} \tag{19}$$

is parameter-free.

Therefore the expected value of the Square Error for the Quantile $\hat{x}_T$

$$QSE(T) = E\left[\left(\hat{x}_T - x_T\right)^2\right] \tag{20}$$

depends on $T$ but not from $a$ and can be calculated from (19) as

$$QSE(T) = b^2 E\left[\left(\hat{z}_T - z_T\right)^2\right] \tag{21}$$

where $z_T$ is the theoretical quantile of the reduced variate $Z$ (6) at the given return period $T$ and $\hat{z}_T$ its graphical (OLS) estimate.

Through a Monte Carlo simulation, $M = 10000$ pseudo-random samples of size $n = 5, 10, 30$ are drawn from Gumbel and Normal parent distributions. Each sample is separately plotted on the corresponding probability paper by using the plotting positions reported in Table 1 and Table 2 and the ones proposed (16).

From (21) we observe that we must compare the $QSE$ only for $b=1$, since it is independent on $a$ (i.e., for each $T$ the statistic $QSE/b^2$ is parameter-free) differently from all the plotting position comparison criteria appeared in the literature: Gringorten (1963), Cunnane (1978), Arnell (1986), Guo (1990), Erto and Lepore (2013). In particular, note that Guo (1990) does not exploit a Monte Carlo simulation, but estimates the quantiles $\hat{x}_T$ on the basis of a single "representative" sample for each plotting position.

In order to compare plotting positions predictive ability independently from $T$, we can consider the integral value of (21), namely Integral Quantile Squared Error ($IQSE$) as follows

$$IQSE = \int_0^1 QSE\left((1-F)^{-1}\right)dF. \tag{22}$$

Moreover, for each parent distribution $F(x;a,b)$ (see Table 3), the predictive ability can be also tested in terms of the expected value of the Squared Error of $F\left(x_T;\hat{a},\hat{b}\right)$

$$FSE(T) = E\left[\left(F\left(x_T;\hat{a},\hat{b}\right) - (1-1/T)\right)^2\right] \tag{23}$$

which depends on $T$ but is parameter-free (Erto 2013). Therefore, as for $QSE$, we can also use its integral value over the $T$ domain

$$IFSE = \int_0^1 FSE\left((1-F)^{-1}\right)dF \tag{24}$$

which is independent from $T$.

**Table 3: Parent distributions**

|  | $F(x;a,b)$ |
|---|---|
| *Gumbel* | $\exp\left[-\exp\left\{-(x-a)/b\right\}\right]$    $0<a, x<+\infty;\ b>0$ |
| *Normal* | $\int_{-\infty}^{x}\left(2\pi b^2\right)^{-1/2}\exp\left\{-(z-a)^2\big/\left(2b^2\right)\right\}dz$ |
| *3-parameter Log-Normal* | $\int_{-\infty}^{x}\left(2\pi b^2\right)^{-1/2}\exp\left\{-\left(\log(x-c)-a\right)^2\big/\left(2b^2\right)\right\}dz$ |

Lastly, it is clear that the measure of plotting position *descriptive ability* introduced by Guo (1990)

$$RM = \sqrt{\frac{1}{N}\sum_{i=1}^{N}\left(\frac{F^{-1}\left(\hat{F}_i;a,b\right)-E\left[X_{(i)}\right]}{E\left[X_{(i)}\right]}\right)^2} \qquad (25)$$

is not parameter-free, and can be modified in light of (18) and (19) into the following index

$$DSE = \sqrt{\frac{1}{N}\sum_{i=1}^{N}\left(\frac{F^{-1}\left(\hat{F}_i;a,b\right)-E\left[X_{(i)}\right]}{b}\right)^2} \qquad (26)$$

which is independent from $a$ and $b$. However, none of the three indices (22), (24) and (26) can be used alone in order to elect the best plotting position, which therefore depends on the purpose for which the results are to be used. Therefore, the average of these three indices for the Gumbel and Normal distributions are reported in Table 4 and confirms the advantages in using the proposed plotting position instead of the classical plotting positions (Table 1 and Table 2). The *IQSE*, *IFSE* and *DSE* values are reported in the Appendix II (Table 9, Table 10 and Table 11, respectively). In particular, the *IQSE* and *IFSE* values are calculated also for the Maximum Likelihood Estimation (MLE) method.

**Table 4: Average of IQSE, IFSE and DSE for the Gumbel and the Normal distributions**

|  | Gumbel | | | Normal | | |
| --- | --- | --- | --- | --- | --- | --- |
|  | $N=5$ | $N=10$ | $N=30$ | $N=5$ | $N=10$ | $N=30$ |
| Erto and Lepore | 0.243 | 0.114 | 0.039 | 0.120 | 0.056 | 0.019 |
| Hazen (1914)-Foster (1936) | 0.260 | 0.129 | 0.048 | 0.139 | 0.072 | 0.026 |
| Beard (1943) | 0.300 | 0.147 | 0.056 | 0.127 | 0.063 | 0.023 |
| Blom (1958) | 0.278 | 0.134 | 0.049 | 0.121 | 0.057 | 0.019 |
| Tukey (1962)-Kerman (2011) | 0.291 | 0.143 | 0.054 | 0.122 | 0.060 | 0.022 |
| Gringorten (1963) | 0.263 | 0.119 | 0.045 | 0.129 | 0.065 | 0.022 |
| Yu and Huang (1999) (a) | 0.300 | 0.149 | 0.057 | 0.142 | 0.072 | 0.026 |
| Yu and Huang (1999) (b) | 0.292 | 0.137 | 0.048 | 0.125 | 0.060 | 0.019 |
| De (2000) | 0.242 | 0.117 | 0.042 | 0.138 | 0.068 | 0.024 |
| Weibull (1914) | 0.432 | 0.217 | 0.088 | 0.193 | 0.101 | 0.041 |
| Cunnane (1977) | 0.271 | 0.130 | 0.047 | 0.124 | 0.060 | 0.020 |
| Adamowski (1981) | 0.323 | 0.160 | 0.063 | 0.139 | 0.071 | 0.027 |
| Erto and Lepore (2013) | 0.303 | 0.148 | 0.057 | 0.129 | 0.064 | 0.023 |

## 4. A critical application

Campi Flegrei is a large volcanic complex located west of the city of Naples, around the town of Pozzuoli (Italy). During the 1983-1984 bradyseismic crisis (slow vertical ground uplift) a total seismic energy of about $4 \cdot 10^{13}$ J (Lima *et al.* 2009) was released. The ground uplift and continuous seismic activity diffused unpleasant emotion and conviction that a volcanic explosion was coming. The "scientific" proof of this upcoming event was given by the Mogi's model (Mogi 1958). This model explains the uplift of a volcanic area as the consequence of the *instability* due to the increasing pressure in the underlying magma that tries to reach the surface. That induced city managers to order a devastating full-scale evacuation of the area. The alternative hypothesis, that explains the ground movement as the consequence of the specific thermo-fluid-dynamics activity of the subsoil of the Campi Flegrei area (Casertano *et al.* 1976), was immediately abandoned. Probably, the careful consideration of

— the time *stability* of the earthquakes magnitude (Figure 1 and Figure 2)

— the complete *independence* of both levels and times of the magnitudes from the focus depths of the corresponding earthquakes (that could have been verified very easily) should have been enough to judge unlikely the hypothesis of an ascending magmatic intrusion, that would have caused ascending rock fractures and consequent ascending earthquake focuses (with time decreasing depths).

**Table 5: Lunar months from July 1983 to July 1984**

|      | I    | II     | III       | IV      | V        | VI       | XIII  |
|------|------|--------|-----------|---------|----------|----------|-------|
| 1983 | July | August | September | October | November | December |       |
|      | 10/07 | 08/08 | 07/09    | 06/10   | 04/11    | 04/12    |       |
|      | 07/08 | 06/09 | 05/10    | 03/11   | 03/12    | 02/01    |       |
|      | VII  | VIII   | IX        | X       | XI       | XII      | XIII  |
| 1984 | January | February | March  | April   | May      | June     | July  |
|      | 03/01 | 02/02 | 02/03    | 01/04   | 01/05    | 30/05    | 29/06 |
|      | 01/02 | 01/03 | 31/03    | 30/04   | 29/05    | 28/06    | 27/07 |

In the Appendix I, the magnitudes (filtered by values less than 1) registered from July 1983 to July 1984 are grouped by lunar month (labelled by I … XIII in Table 5) because of the high correlation among bradyseism and short and long period tidal components (Casertano *et al.* 1976). In Figure 1, the data are analysed via 3-parameter Log-Normal (Table 3) probability paper by using the proposed plotting positions (16) with $k = 4$ and threshold parameter $c = 1$ (Table 3). Some of the step increases of $y_{(i)} = E\{Z_{(i)}\}$ are wider because of the ties resulting from the measurement resolution. This real scenario is one of the typical critical cases where a reliable graphical analysis of the data is the only persuasive way to share statistical conclusions with non-statistician managers that have to utilize them to make grave decisions on territory and citizens. Table 6 reports from Figure 1 the probability of a magnitude greater than 5, which in expert opinion is the critical threshold for concrete structures. Figure 1 and Figure 2 are evident pictures of substantial earthquakes *stability* during the 1983-1984 bradyseismic crisis. Obviously, they could have helped to warn against the alarmism caused by the apocalyptical newspaper titles at the time (Gore and Mazzatenta 1984). No analytical results or indices could have provided so rapid information to non-statisticians.

However, all the graphical results are confirmed by the analytical goodness-of-fit tests for (log)normality carried out also through the modified Anderson-Darling (mAD) upper-tail test (Stephens1974, Case 3). In particular, Table 7 reports the mAD statistic values to test the goodness-of-fit of the (Log)Normal distribution for the magnitudes greater than 1 for each lunar month whereas Table 8 reports the mAD statistic values to test whether such magnitudes belongs also to the (Log)Normal distribution with the population (unknown) parameters estimated on the basis of the cumulative sample of the previous months. Most the values reported in these tables do not exceed the significance point at level 5.0 (0.787). In any case, the maximum value (0.825) does not exceed the significance point at level 2.5 (0.918).

**Table 6: Probability estimates of a magnitude greater than 5**

| Lunar month | | Lunar month | |
|---|---|---|---|
| I | 0.0001 | VII | 0.0020 |
| II | 0.0017 | VIII | 0.0013 |
| III | 0.0060 | IX | 0.0037 |
| IV | 0.0012 | X | 0.0021 |
| V | 0.0028 | XI | 0.0013 |
| VI | 0.0060 | XII | 0.0058 |
| | | XIII | 0.0018 |

**Table 7: Modified Anderson-Darling statistic values with population (unknown) parameters estimated at each lunar month**

| Lunar month | | Lunar month | |
|:---:|:---:|:---:|:---:|
| **I** | 0.521 | **VII** | 0.701 |
| **II** | 0.313 | **VIII** | 0.732 |
| **III** | 0.766 | **IX** | 0.466 |
| **IV** | 0.382 | **X** | 0.709 |
| **V** | 0.577 | **XI** | 0.754 |
| **VI** | 0.424 | **XII** | 0.434 |
| | | **XIII** | 0.494 |

**Table 8: Modified Anderson-Darling statistic values with population (unknown) parameters estimated on the basis of the cumulative sample of the previous lunar months**

| Lunar month | | Lunar month | |
|:---:|:---:|:---:|:---:|
| **I** | 0.521 | **VII** | 0.752 |
| **II** | 0.395 | **VIII** | 0.825 |
| **III** | 0.744 | **IX** | 0.821 |
| **IV** | 0.635 | **X** | 0.739 |
| **V** | 0.710 | **XI** | 0.742 |
| **VI** | 0.702 | **XII** | 0.699 |
| | | **XIII** | 0.680 |

**Figure 1: Substantial stability of the distribution of magnitudes from July 1983 (I) to July 1984 (XIII) (see Table 5) graphically shown**

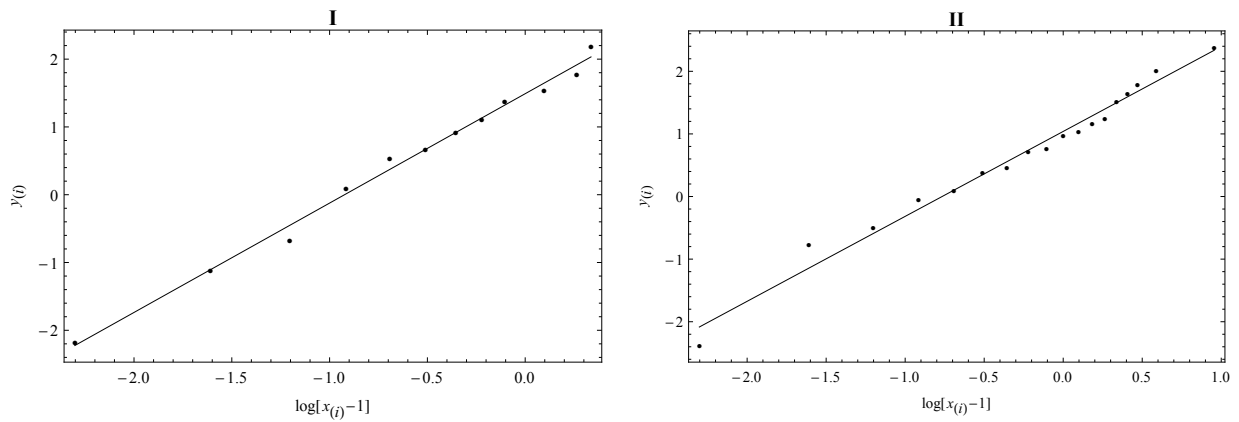

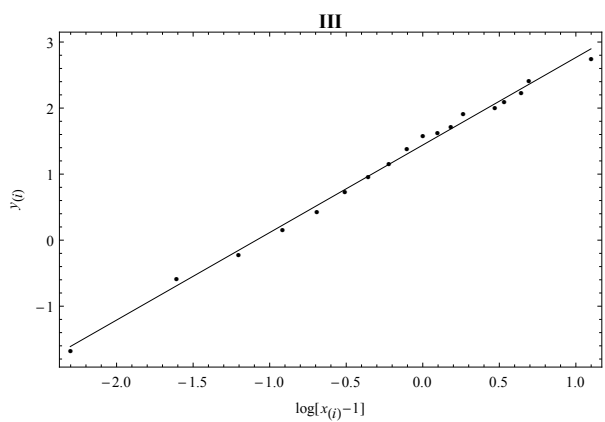

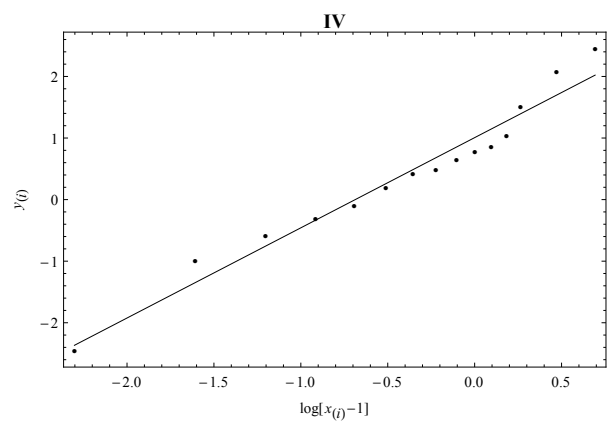

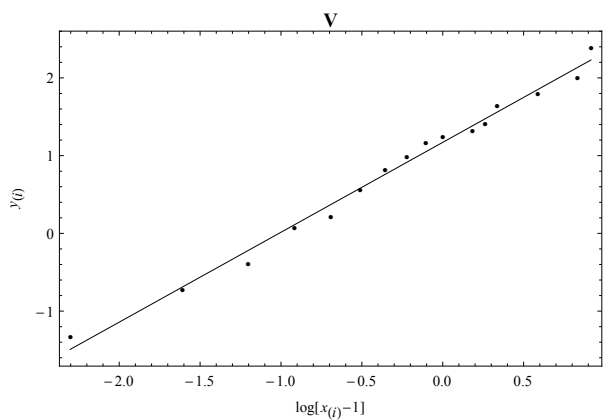

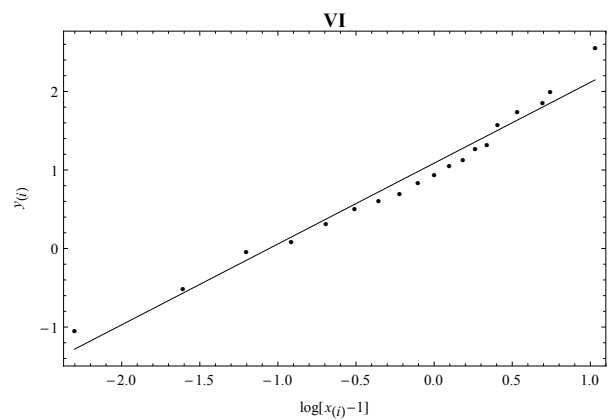

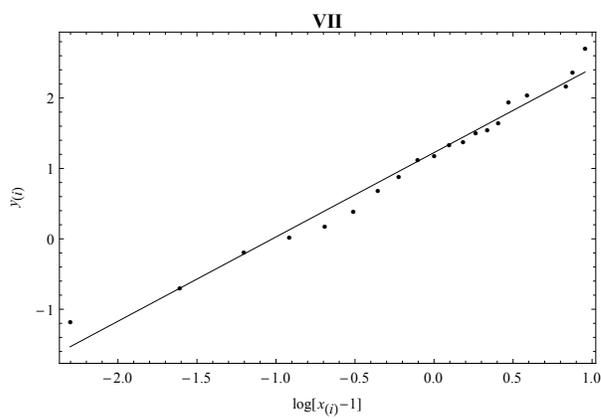

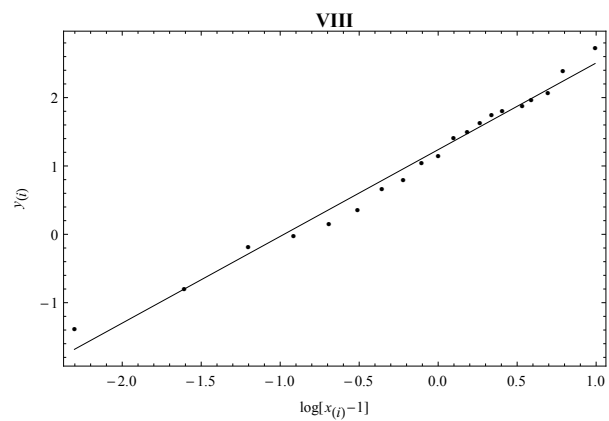

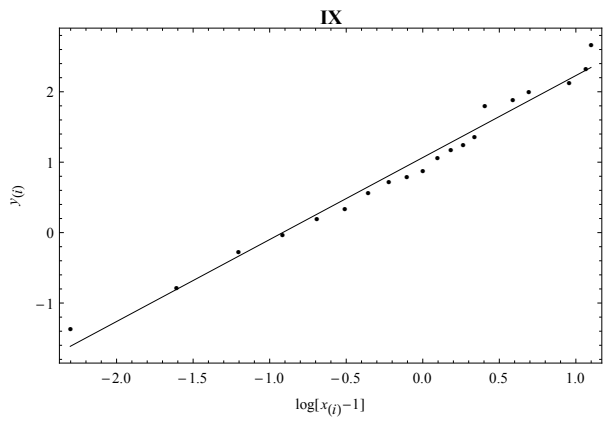

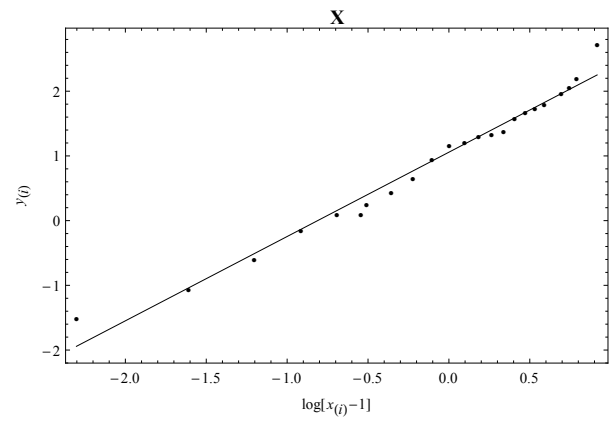

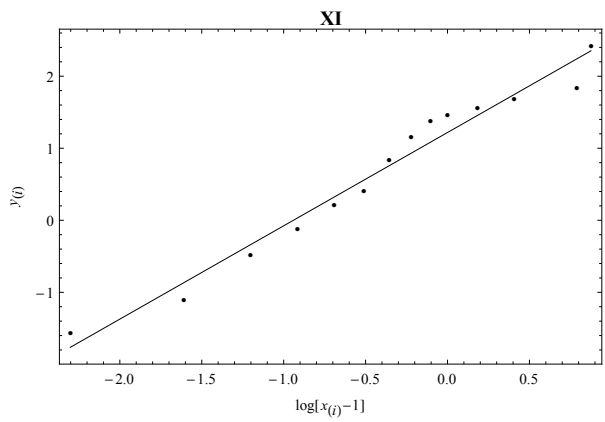

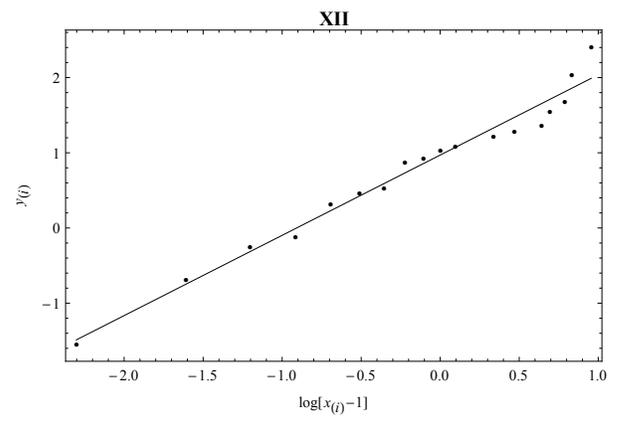

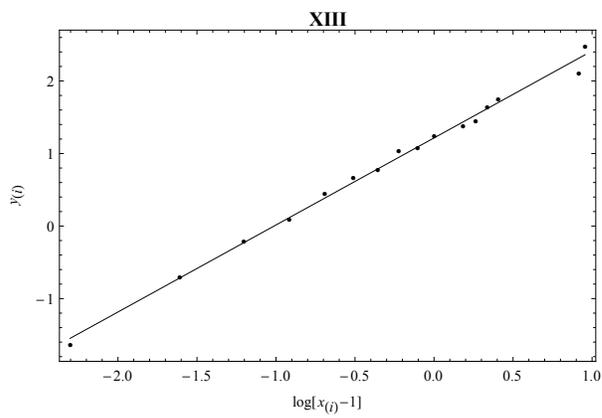

**Figure 2: Picture (from Figure1) of the non-increasing probability of a magnitude greater than 5 from July 1983 (I) to July 1984 (XIII)**

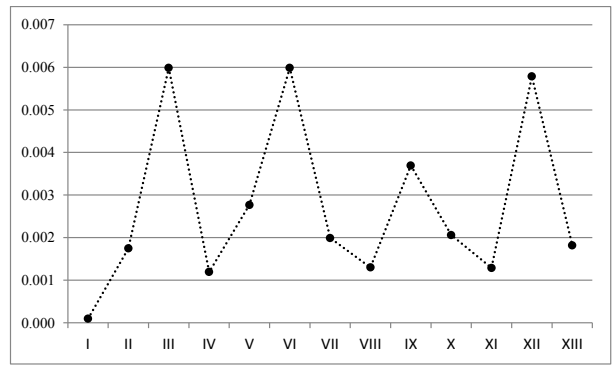

## 5. Conclusions

On the basis of theoretical considerations, a new plotting position approach is proposed and is compared with classical ones (Table 1 and Table 2). In order to test plotting position descriptive and predictive ability, new parameter-free performance indices ( $DSE$ , $IQSE$ and $IFSE$ ) are adopted. These indices simplify and generalize those known in the literature. A wide Monte Carlo simulation confirms that the proposed plotting position formula, which leads to BLUE of location and scale parameters, out-performs all the classical plotting positions (Table 4) even by considering only the first $k = 4$ terms of the Taylor expanding (9) on which it is based. As $k$ increases, it can reach any precision in terms of $DSE$ (Table 11). Moreover, the proposed formula shows $IQSE$ and $IFSE$ values (that test plotting position predictive ability) which are always very close to the *best values* (Table 9 and Table 10) differently from those shown by all the classical plotting positions (Table 1 and Table 2). Note that the benchmarking *best values* come from a formula always different from case to case. The results are obtained both in the case of skewed (Gumbel) and simmetrycal (Normal) parent distributions. The good properties of the proposed plotting position formula reduce the efficiency/reputation gap between probability plotting and the corresponding analytical methods. That encourages adopting graphical procedures, based on probability plotting, which result very useful in critical application where visual representation of the results of statistical analyses help non-statisticians to make right decisions.

## Appendix I:  the bradyseism magnitudes (filtered by values less than 1) registered during in Campi Flegrei (Italy) from July 1983 to July 1984.

*July 1983*

1.3, 1.5, 1.1, 1.2, 2.3, 1.4, 1.2, 1.2, 1.4, 1.4, 1.4, 1.0, 1.4, 1.0, 1.2, 1.2, 1.4, 1.6, 1.7, 1.3, 1.5, 1.5, 1.7, 1.3, 1.5, 1.4, 1.4, 1.0, 1.5, 1.3, 2.4, 1.6, 1.8, 1.4, 1.9, 2.1, 1.3, 1.8, 1.7, 1.4, 1.0, 1.9, 1.5, 1.5, 1.4, 1.4

*August 1983*

1.4, 2.0, 2.2, 1.0, 1.8, 2.4, 1.1, 2.3, 1.5, 1.0, 1.2, 1.3, 1.2, 2.1, 1.2, 1.3, 1.4, 1.6, 1.5, 1.6, 1.7, 2.8, 1.7, 1.5, 1.8, 1.2, 1.4, 1.0, 1.8, 1.4, 1.6, 1.6, 1.0, 1.4, 1.2, 1.2, 1.2, 1.2, 2.0, 1.9, 1.6, 1.5, 1.3, 1.6, 1.4, 2.5, 2.4, 1.6, 2.0, 1.4, 1.8, 1.3, 1.3, 1.2, 1.8, 1.2, 1.2, 1.2, 1.2, 1.2, 2.2, 1.0, 1.4, 1.4, 2.4, 1.4, 3.6, 1.8, 1.3, 1.6, 2.0, 1.4, 1.0, 1.4, 2.6, 1.2

*September 1983*

2.0, 1.6, 1.4, 1.0, 1.4, 1.3, 1.5, 1.7, 1.2, 1.0, 1.4, 1.0, 1.1, 1.4, 1.0, 1.4, 1.2, 1.0, 1.3, 1.6, 1.2, 1.6, 1.2, 1.2, 1.2, 2.0, 1.0, 1.2, 1.2, 1.5, 1.0, 1.8, 1.1, 1.2, 1.4, 1.6, 1.8, 1.2, 1.9, 1.4, 1.2, 1.4, 1.4, 2.2, 1.2, 2.6, 1.0, 1.5, 1.2, 1.8, 2.7, 1.8, 1.2, 1.0, 1.5, 1.2, 1.4, 1.4, 1.4, 1.5, 1.0, 1.6, 1.5, 1.2, 1.5, 1.4, 1.7, 1.2, 1.5, 1.2, 1.7, 1.2, 1.2, 1.3, 1.4, 1.3, 2.3, 1.7, 1.8, 1.2, 1.4, 1.3, 1.4, 2.0, 1.3, 1.6, 1.6, 1.4, 1.2, 1.6, 1.3, 1.4, 1.0, 1.2, 1.0, 1.1, 1.1, 1.8, 1.5, 1.2, 1.9, 1.0, 1.7, 1.2, 1.0, 1.2, 1.0, 1.3, 1.2, 1.5, 2.3, 1.0, 1.2, 1.4, 1.3, 1.0, 1.2, 1.4, 1.0, 1.4, 1.2, 1.0, 1.0, 1.8, 1.5, 1.5, 1.1, 1.6, 1.1, 1.0, 1.0, 1.9, 1.0, 1.2, 1.5, 1.2, 1.1, 1.2, 1.1, 1.9, 1.3, 1.2, 1.9, 1.5, 1.0, 1.0, 1.3, 1.4, 1.0, 1.2, 1.4, 1.5, 1.2, 1.1, 1.7, 1.1, 1.4, 1.2, 1.2, 1.9, 1.5, 1.0, 1.2, 1.2, 1.7, 1.0, 1.6, 1.0, 1.3, 1.4, 2.0, 1.7, 2.3, 1.3, 2.9, 1.7, 1.8, 1.6, 1.7, 1.6, 1.2, 1.6, 1.6, 1.0, 1.6, 2.0, 1.0, 1.7, 1.0, 1.5, 1.4, 1.8, 1.0, 1.3, 1.5, 1.6, 1.5, 1.4, 1.7, 1.5, 1.6, 1.3, 1.0, 1.0, 1.3, 1.5, 1.4, 1.4, 1.3, 1.0, 1.6, 1.7, 1.6, 1.2, 1.2, 1.2, 1.3, 1.9, 2.1, 1.2, 1.3, 1.3, 1.3, 2.2, 1.5, 1.3, 1.3, 1.9, 1.0, 1.4, 1.2, 4.0, 1.2, 1.3, 1.3, 1.6, 1.3, 1.0, 3.0

*October 1983*

1.5, 1.2, 1.9, 1.4, 2.2, 1.0, 1.5, 1.6, 1.3, 1.3, 1.3, 1.0, 1.2, 2.0, 1.0, 1.2, 1.0, 2.3, 1.9, 1.3, 1.5, 1.5, 1.6, 1.3, 2.0, 1.0, 1.5, 1.0, 2.3, 1.2, 1.5, 1.7, 1.8, 1.2, 2.0, 1.0, 1.4, 1.3, 1.3, 1.4, 1.1, 1.4, 1.6, 1.0, 3.0, 2.1, 2.6, 2.3, 1.2, 2.3, 1.0, 2.3, 1.9, 1.6, 2.6, 2.6, 1.0, 2.2, 1.4, 1.2, 1.0, 1.4, 2.2, 1.3, 1.7, 1.9, 1.9, 2.3, 1.4, 1.6, 1.7, 1.2, 1.2, 1.5, 1.0, 1.6, 2.1, 2.2, 1.0, 1.6, 1.5, 1.7, 1.7, 1.4, 1.6, 1.6, 1.0, 2.3, 1.3, 1.6, 1.2, 1.7, 1.2, 1.7, 1.2, 1.8, 1.0, 1.3, 1.2, 1.0, 2.6

*November 1983*

1.4, 1.3, 2.0, 1.1, 1.2, 1.7, 1.9, 1.7, 1.5, 1.8, 1.6, 1.1, 1.1, 1.1, 1.1, 1.6, 1.4, 1.2, 1.2, 1.3, 1.7, 1.8,

1.4, 2.2, 1.6, 1.6, 1.2, 1.4, 1.8, 1.3, 1.5, 1.7, 2.8, 1.6, 1.4, 1.5, 1.4, 3.3, 1.4, 1.3, 1.4, 1.2, 1.6, 1.1, 1.0, 1.7, 1.5, 1.2, 1.2, 1.3, 2.4, 1.4, 1.0, 1.3, 1.3, 1.2, 1.6, 1.0, 1.2, 3.5, 1.0, 1.4, 1.6, 1.4, 1.2, 1.9, 1.9, 1.0, 1.1, 1.0, 1.6, 1.3, 1.7, 2.4, 1.4, 2.3, 1.0, 1.4, 1.0, 1.0

### December 1983

1.0, 1.9, 1.0, 1.1, 1.0, 1.1, 1.5, 1.5, 1.7, 1.3, 1.1, 1.1, 1.2, 1.4, 1.5, 1.0, 1.2, 2.1, 1.3, 1.2, 2.0, 1.9, 1.2, 1.7, 1.8, 1.8, 1.2, 1.5, 1.3, 1.1, 1.1, 1.5, 1.9, 1.2, 2.2, 1.4, 1.4, 1.5, 1.8, 1.4, 1.1, 2.1, 1.3, 1.3, 1.1, 1.0, 1.1, 1.2, 1.0, 1.0, 2.4, 2.1, 2.5, 2.7, 1.2, 1.3, 1.4, 3.1, 1.2, 2.3, 1.6, 1.6, 1.2, 1.2, 1.6, 1.0, 1.5, 3.8, 1.2, 1.9, 1.5, 1.1, 1.3, 1.3, 1.0, 1.3, 2.5, 3.0, 1.0, 1.2, 1.7, 2.5, 1.3, 1.1, 1.3, 1.0, 1.1, 1.0, 1.3, 1.2, 1.6, 1.0, 1.0, 1.5, 1.0, 1.3, 1.2, 1.1, 1.4, 1.3, 1.6, 1.9, 2.2, 2.7, 3.8, 1.7, 1.6, 2.3, 1.1, 1.1, 1.2, 1.3, 2.3, 2.0, 1.6, 1.3, 2.0, 1.3, 1.5, 1.2, 1.6, 1.0, 1.2, 2.5, 1.0, 1.3, 1.1, 1.3, 1.3, 1.1

### January 1984

1.1, 1.5, 1.2, 1.7, 2.5, 1.5, 1.7, 1.5, 1.3, 2.1, 1.9, 1.1, 1.0, 1.9, 1.6, 1.3, 1.3, 1.0, 1.0, 1.1, 1.3, 1.9, 1.3, 1.7, 1.2, 1.1, 2.3, 1.3, 1.3, 1.1, 2.0, 1.4, 1.3, 1.1, 1.1, 1.6, 1.1, 1.1, 1.1, 1.2, 1.8, 1.1, 1.3, 1.3, 1.2, 2.6, 1.6, 1.0, 1.3, 2.4, 1.0, 1.0, 1.4, 1.5, 1.8, 1.6, 1.4, 1.8, 1.2, 3.3, 1.2, 1.2, 1.1, 1.7, 2.0, 1.3, 1.5, 1.4, 1.6, 1.7, 1.3, 1.1, 1.1, 1.8, 1.0, 1.4, 1.2, 1.1, 1.2, 1.0, 1.1, 1.1, 1.0, 1.2, 1.4, 1.3, 1.5, 1.7, 2.8, 2.6, 1.6, 1.3, 1.4, 1.0, 1.3, 1.6, 1.3, 1.8, 1.3, 1.7, 1.9, 1.0, 1.5, 1.3, 1.2, 1.3, 1.2, 1.8, 1.2, 2.1, 1.9, 1.4, 1.5, 1.7, 2.6, 1.1, 1.3, 1.4, 1.8, 2.1, 1.3, 1.3, 1.1, 1.4, 1.5, 1.5, 1.4, 3.4, 1.9, 1.4, 1.8, 1.3, 2.1, 1.6, 2.2, 1.9, 1.2, 2.3, 1.7, 1.6, 1.6, 1.2, 1.2, 1.6, 1.7, 1.3, 2.6, 1.7, 1.0, 1.9, 1.3, 1.4, 1.1, 1.3, 1.7, 1.3, 1.9, 1.1, 1.0, 2.5, 1.7, 1.6, 1.2, 1.7, 1.7, 1.7, 1.2, 1.5, 1.6, 1.8, 1.0, 1.3, 2.3, 1.3, 1.7, 1.6, 1.3, 2.1, 1.2, 1.7, 1.8, 3.6, 1.9, 1.4, 1.2, 1.3, 1.2

### February 1984

1.1, 1.4, 1.8, 1.3, 1.0, 1.6, 1.6, 1.5, 1.4, 1.9, 1.7, 1.5, 1.0, 1.1, 1.4, 1.4, 1.0, 1.2, 1.1, 1.1, 1.5, 1.3, 1.7, 1.3, 2.1, 1.4, 1.3, 1.7, 1.3, 1.4, 1.0, 1.3, 1.6, 1.0, 2.4, 1.3, 1.4, 1.3, 1.3, 1.3, 1.2, 1.3, 1.1, 1.3, 1.0, 1.2, 1.3, 1.3, 1.5, 1.9, 1.3, 1.2, 1.8, 1.6, 1.5, 1.2, 1.3, 1.2, 1.5, 1.3, 1.1, 1.3, 2.1, 1.3, 3.2, 1.9, 1.2, 2.1, 1.6, 1.6, 1.8, 1.8, 2.7, 2.5, 2.1, 1.7, 2.3, 2.1, 1.9, 2.1, 1.8, 1.0, 1.1, 1.9, 1.3, 1.7, 1.6, 1.7, 2.0, 2.3, 1.6, 1.2, 1.0, 1.3, 1.6, 1.9, 1.9, 1.3, 1.8, 1.2, 1.7, 1.1, 1.0, 1.9, 1.3, 1.5, 1.5, 1.3, 1.3, 1.2, 1.9, 1.1, 1.9, 1.7, 1.6, 1.7, 1.2, 1.7, 1.8, 1.4, 1.6, 1.4, 2.4, 1.9, 1.6, 1.7, 1.4, 1.3, 2.8, 1.6, 1.5, 1.7, 1.3, 1.2, 1.2, 2.0, 1.7, 3.7, 1.3, 1.5, 1.7, 1.0, 1.1, 1.2, 1.3, 1.3, 2.2, 1.9, 1.3, 1.2, 1.5, 1.3, 1.4, 1.6, 2.1, 1.7, 1.0, 1.3, 1.3, 1.0, 3.0, 3.2, 1.0, 1.0, 1.0, 1.0, 1.2, 1.1, 1.2, 1.2, 1.4, 1.3, 2.3, 1.7, 2.0, 1.3, 2.1, 1.1, 1.2, 1.7, 1.6, 1.3, 1.1, 1.7, 2.1, 1.3, 1.7, 1.8, 2.2, 1.0, 2.0, 1.2, 1.2, 1.0, 1.7, 1.2, 1.2, 1.0, 1.5, 1.3, 1.1, 1.3

### March 1984

1.8, 2.5, 1.4, 1.0, 1.3, 1.2, 1.3, 1.7, 1.5, 1.0, 1.2, 1.5, 2.5, 1.2, 1.2, 1.2, 2.5, 1.7, 1.8, 1.4, 1.2, 1.5, 1.2, 2.2, 2.1, 1.6, 2.3, 1.4, 1.1, 1.1, 1.3, 1.7, 1.5, 1.6, 1.0, 1.2, 1.0, 1.3, 1.8, 1.5, 1.7, 1.3, 1.8, 1.5, 1.4, 2.2, 1.3, 1.5, 1.1, 2.1, 1.0, 3.9, 1.4, 1.3, 1.4, 1.1, 1.0, 1.1, 1.2, 1.0, 1.8, 1.0, 2.0, 1.0, 1.0, 1.4, 1.0, 1.7, 1.3, 1.2, 1.2, 1.8, 1.0, 2.5, 1.2, 1.5, 1.4, 1.7, 1.2, 2.8, 1.9, 1.6, 2.5, 1.0, 1.8, 2.0, 1.6, 1.0, 1.5, 1.3, 1.3, 2.5, 1.3, 1.9, 1.8, 1.7, 2.1, 4.0, 1.1, 1.3, 1.5, 2.1, 1.0, 1.1, 1.0, 2.4, 2.0, 1.4, 1.6, 1.1, 1.3, 1.2, 1.3, 1.2, 1.3, 3.6, 2.5, 2.2, 1.3, 1.1, 2.1, 1.1, 1.4, 1.3, 1.2, 1.7, 1.0, 1.7, 1.1, 1.3, 1.3, 1.5, 1.4, 2.4, 1.4, 1.0, 1.3, 1.7, 1.2, 2.5, 3.0, 2.4, 1.6, 1.7, 2.1, 1.0, 1.5, 1.5, 1.0, 1.4, 1.3, 1.9, 1.1, 1.0, 2.3, 1.0, 1.1, 1.2, 1.3, 1.4, 2.0, 1.4, 1.3, 1.6, 1.2, 1.3, 1.1, 1.7, 1.5, 2.2, 1.3, 1.7, 1.0, 1.3, 1.3, 1.3, 2.1, 1.6, 1.2

### April 1984

1.0, 1.7, 1.0, 1.4, 1.6, 1.3, 1.8, 1.9, 1.3, 1.8, 1.4, 1.0, 1.5, 1.9, 1.4, 1.9, 2.3, 1.5, 1.3, 1.4, 1.9, 1.5, 1.5, 1.6, 2.0, 2.0, 1.4, 2.2, 1.4, 1.8, 1.6, 1.4, 2.0, 1.0, 1.7, 1.8, 2.7, 1.3, 2.5, 1.6, 3.0, 1.4, 1.4, 1.8, 1.8, 1.7, 1.2, 1.7, 1.8, 3.0, 1.1, 1.9, 1.0, 1.8, 2.5, 2.5, 1.4, 1.3, 1.3, 1.9, 1.4, 1.5, 1.5, 1.7, 1.4, 2.5,

2.0, 1.1, 1.5, 1.8, 1.2, 2.2, 1.6, 1.6, 1.3, 1.1, 1.1, 1.1, 1.0, 2.1, 1.9, 1.8, 1.3, 1.8, 1.5, 1.5, 1.4, 1.4, 1.5, 1.0, 1.5, 1.2, 1.9, 1.0, 1.7, 1.1, 1.7, 1.2, 1.1, 1.2, 1.0, 1.4, 1.8, 1.8, 1.0, 1.1, 1.6, 1.4, 1.0, 1.0, 1.0, 1.0, 1.4, 1.5, 1.3, 1.3, 1.5, 1.7, 1.2, 3.5, 1.3, 1.4, 1.3, 1.3, 2.0, 1.7, 1.4, 1.2, 1.3, 1.3, 2.0, 2.0, 1.2, 1.3, 1.7, 1.3, 2.6, 2.0, 1.2, 1.5, 1.3, 1.4, 1.5, 1.0, 1.1, 1.9, 1.6, 1.9, 1.9, 1.0, 1.7, 1.0, 1.3, 1.5, 2.6, 1.9, 1.4, 1.9, 1.0, 1.0, 1.0, 1.0, 1.1, 1.0, 2.0, 1.4, 1.0, 1.9, 1.0, 1.4, 1.1, 1.0, 1.4, 1.4, 1.0, 1.9, 1.8, 1.3, 1.0, 1.3, 2.8, 1.2, 1.0, 1.5, 1.3, 2.5, 1.6, 1.3, 3.5, 1.4, 1.4, 1.4, 1.0, 1.1, 1.5, 1.2, 1.2, 1.6, 1.7, 1.4, 3.1, 2.4, 3.2, 1.2, 1.7, 1.2, 2.1, 2.2, 1.0, 1.4, 1.3

**May 1984**

1.7, 1.4, 1.0, 1.4, 1.5, 1.9, 1.2, 1.4, 1.0, 1.8, 1.7, 1.0, 1.3, 1.9, 1.0, 1.5, 1.3, 1.6, 1.9, 1.0, 3.4, 1.2, 1.0, 2.5, 1.7, 1.8, 1.4, 3.4, 1.3, 1.4, 1.1, 1.1, 1.0, 1.3, 1.0, 1.3, 1.3, 1.3, 1.0, 1.5, 1.7, 1.7, 1.4, 1.3, 1.5, 1.1, 1.3, 3.2, 1.2, 1.7, 1.6, 1.7, 1.4, 1.3, 1.0, 1.4, 1.7, 1.3, 1.0, 2.2, 1.1, 1.5, 1.6, 2.0, 1.4, 1.2, 1.1, 1.3, 1.5, 1.7, 1.4, 1.0, 1.7, 1.4, 1.8, 1.8, 1.5, 1.5, 1.6, 1.5, 1.6, 1.6, 1.3, 1.0, 1.2, 1.2, 1.3, 1.5, 1.7, 1.8, 1.8

**June 1984**

1.5, 3.2, 1.5, 1.3, 3.3, 1.8, 3.0, 1.5, 1.3, 1.8, 1.8, 1.5, 3.0, 1.2, 1.3, 1.2, 1.5, 2.0, 2.1, 1.8, 1.2, 1.8, 1.2, 1.5, 1.3, 1.3, 1.5, 1.7, 1.6, 2.0, 1.3, 1.6, 1.0, 1.4, 1.8, 1.2, 1.8, 1.5, 1.2, 1.5, 1.9, 1.1, 1.5, 1.2, 1.3, 1.3, 1.0, 1.3, 2.4, 2.4, 2.9, 1.4, 1.4, 1.3, 1.5, 1.0, 1.0, 1.1, 1.2, 1.4, 1.6, 1.1, 1.0, 1.0, 1.2, 1.2, 1.0, 1.5, 1.2, 1.5, 1.1, 1.2, 1.0, 1.0, 2.6, 3.3, 1.2, 1.6, 1.3, 1.0, 1.2, 1.3, 1.1, 1.7, 1.0, 3.6, 1.8

**July 1984**

1.0, 1.5, 1.3, 1.4, 1.3, 1.3, 1.1, 1.5, 1.5, 1.5, 1.0, 1.3, 1.5, 3.6, 2.2, 1.0, 1.5, 1.4, 1.9, 3.5, 1.6, 1.2, 1.8, 1.2, 1.6, 1.7, 1.4, 1.6, 1.5, 1.8, 1.6, 1.4, 1.2, 1.5, 2.4, 1.6, 1.0, 1.3, 1.3, 1.2, 1.2, 1.2, 1.0, 1.0, 1.0, 1.3, 1.5, 1.3, 1.0, 1.3, 1.5, 1.0, 1.2, 1.6, 1.0, 1.0, 3.5, 1.0, 1.3, 1.3, 2.0, 1.2, 1.4, 1.4, 1.2, 1.0, 1.8, 1.3, 1.5, 1.0, 1.4, 1.4, 1.2, 1.2, 1.2, 2.0, 2.3, 1.7, 1.2, 1.0, 1.0, 1.2, 1.4, 1.8, 1.0, 2.4, 1.7, 1.3, 2.2, 1.1, 1.3, 1.0, 1.4, 1.1, 1.0, 2.0, 1.6, 2.5, 1.4, 1.0, 1.0, 1.1, 1.2, 1.5, 1.3, 1.2, 1.2, 1.8, 1.3, 1.1

## Appendix II: The *IQSE*, *IFSE* and *DSE* for the Gumbel and Normal distributions

**Table 9: IQSE for the Gumbel and Normal distributions**

|  | Gumbel | | | Normal | | |
|---|---|---|---|---|---|---|
|  | $N = 5$ | $N = 10$ | $N = 30$ | $N = 5$ | $N = 10$ | $N = 30$ |
| MLE | 0.575 | 0.270 | 0.092 | 0.317 | 0.154 | 0.051 |
| Erto and Lepore | 0.693 | 0.326 | 0.113 | 0.334 | 0.157 | 0.051 |
| Hazen (1914)-Foster (1936) | 0.670 | 0.321 | 0.113 | 0.317 | 0.154 | 0.051 |
| Beard (1943) | 0.770 | 0.353 | 0.118 | 0.341 | 0.160 | 0.051 |
| Blom (1958) | 0.730 | 0.339 | 0.116 | 0.330 | 0.157 | 0.051 |
| Tukey (1962) - Kerman (2011) | 0.755 | 0.348 | 0.117 | 0.336 | 0.159 | 0.051 |
| Gringorten (1963) | 0.696 | 0.328 | 0.114 | 0.322 | 0.155 | 0.051 |
| Yu and Huang (1999) (a) | 0.717 | 0.335 | 0.115 | 0.326 | 0.156 | 0.051 |
| Yu and Huang (1999) (b) | 0.777 | 0.353 | 0.118 | 0.329 | 0.158 | 0.051 |
| De (2000) | 0.690 | 0.328 | 0.114 | 0.334 | 0.158 | 0.051 |
| Weibull (1914) | 1.039 | 0.448 | 0.137 | 0.430 | 0.189 | 0.057 |
| Cunnane (1977) | 0.716 | 0.335 | 0.115 | 0.326 | 0.156 | 0.051 |
| Adamowski (1981) | 0.813 | 0.368 | 0.121 | 0.353 | 0.164 | 0.052 |
| Erto and Lepore (2013) | 0.776 | 0.354 | 0.119 | 0.342 | 0.160 | 0.051 |

**Table 10: IFSE for the Gumbel and Normal distributions**

|  | Gumbel | | | Normal | | |
|---|---|---|---|---|---|---|
|  | $N = 5$ | $N = 10$ | $N = 30$ | $N = 5$ | $N = 10$ | $N = 30$ |
| MLE | 0.027 | 0.012 | 0.004 | 0.025 | 0.012 | 0.004 |
| Erto and Lepore | 0.024 | 0.012 | 0.004 | 0.022 | 0.011 | 0.004 |
| Hazen (1914) - Foster (1936) | 0.026 | 0.012 | 0.004 | 0.023 | 0.011 | 0.004 |
| Beard (1943) | 0.024 | 0.011 | 0.004 | 0.021 | 0.011 | 0.004 |
| Blom (1958) | 0.024 | 0.012 | 0.004 | 0.022 | 0.011 | 0.004 |
| Tukey (1962) - Kerman (2011) | 0.024 | 0.012 | 0.004 | 0.022 | 0.011 | 0.004 |
| Gringorten (1963) | 0.025 | 0.012 | 0.004 | 0.022 | 0.012 | 0.004 |
| Yu and Huang (1999) (a) | 0.025 | 0.012 | 0.004 | 0.022 | 0.011 | 0.004 |
| Yu and Huang (1999) (b) | 0.025 | 0.012 | 0.004 | 0.023 | 0.011 | 0.004 |
| De (2000) | 0.024 | 0.012 | 0.004 | 0.022 | 0.011 | 0.004 |
| Weibull (1914) | 0.022 | 0.011 | 0.004 | 0.020 | 0.011 | 0.003 |
| Cunnane (1977) | 0.025 | 0.012 | 0.004 | 0.022 | 0.011 | 0.004 |
| Adamowski (1981) | 0.023 | 0.011 | 0.004 | 0.021 | 0.011 | 0.004 |
| Erto and Lepore (2013) | 0.024 | 0.011 | 0.004 | 0.021 | 0.011 | 0.004 |

**Table 11: DSE for the Gumbel and Normal distributions**

|  | Gumbel | | | Normal | | |
|---|---|---|---|---|---|---|
|  | $N = 5$ | $N = 10$ | $N = 30$ | $N = 5$ | $N = 10$ | $N = 30$ |
| Erto and Lepore | 0.003 | 0.000 | 0.001 | 0.011 | 0.004 | 0.001 |
| Hazen (1914) - Foster (1936) | 0.077 | 0.051 | 0.024 | 0.085 | 0.054 | 0.027 |
| Beard (1943) | 0.019 | 0.019 | 0.014 | 0.105 | 0.077 | 0.046 |
| Blom (1958) | 0.011 | 0.004 | 0.002 | 0.079 | 0.052 | 0.028 |
| Tukey (1962) - Kerman (2011) | 0.009 | 0.011 | 0.010 | 0.095 | 0.068 | 0.040 |
| Gringorten (1963) | 0.043 | 0.027 | 0.011 | 0.068 | 0.016 | 0.016 |
| Yu and Huang (1999) (a) | 0.079 | 0.049 | 0.023 | 0.158 | 0.101 | 0.052 |
| Yu and Huang (1999) (b) | 0.022 | 0.012 | 0.003 | 0.073 | 0.045 | 0.022 |
| De (2000) | 0.058 | 0.036 | 0.018 | 0.012 | 0.011 | 0.008 |
| Weibull (1914) | 0.130 | 0.103 | 0.062 | 0.236 | 0.192 | 0.123 |
| Cunnane (1977) | 0.023 | 0.012 | 0.004 | 0.072 | 0.044 | 0.022 |
| Adamowski (1981) | 0.044 | 0.037 | 0.024 | 0.132 | 0.102 | 0.063 |
| Erto and Lepore (2013) | 0.023 | 0.021 | 0.014 | 0.109 | 0.080 | 0.048 |